\pdfoutput=1

\documentclass[10pt, aps, reprint, numerical, showkeys, showpacs, superscriptaddress, amsmath, amssymb, floatfix]{revtex4-1}

\usepackage[utf8]{inputenc}
\usepackage{amsmath}
\usepackage{amssymb}
\usepackage{physics}
\usepackage{amsthm} 
\usepackage{lmodern}
\usepackage{xcolor}
\usepackage[pdftex]{graphicx}
\usepackage[hidelinks]{hyperref}
\usepackage{csquotes} 
\usepackage{float} 
\usepackage{cleveref} 
\usepackage{verbatim}
\usepackage{enumitem}

\begin{document}

\newcommand\at[2]{\left.#1\right|_{#2}}
\newcommand{\hz}{ \ \text{s}^{-1}}
\newcommand{\hzz}{ \ \text{s}^{-2}}

\title{Delay master stability of inertial oscillator networks}

\author{Reyk B\"orner}
\affiliation{Potsdam Institute for Climate Impact Research (PIK), Member of the Leibniz Association, P.O. Box 60 12 03, D-14412 Potsdam, Germany}
\affiliation{Department of Physics, Freie Universität Berlin, Arnimallee 14, 14195 Berlin, Germany}

\author{Paul \surname{Schultz}}
\affiliation{Potsdam Institute for Climate Impact Research (PIK), Member of the Leibniz Association, P.O. Box 60 12 03, D-14412 Potsdam, Germany}

\author{Benjamin \"Unzelmann}
\affiliation{Department of Physics, Freie Universität Berlin, Arnimallee 14, 14195 Berlin, Germany}
\affiliation{Potsdam Institute for Climate Impact Research (PIK), Member of the Leibniz Association, P.O. Box 60 12 03, D-14412 Potsdam, Germany}

\author{Deli Wang}
\affiliation{School of Science, Xi'an University of Architecture and Technology, Xi’an 710055, China}

\author{Frank Hellmann}\email[Email: ]{hellmann@pik-potsdam.de}
\affiliation{Potsdam Institute for Climate Impact Research (PIK), Member of the Leibniz Association, P.O. Box 60 12 03, D-14412 Potsdam, Germany}

\author{J\"{u}rgen \surname{Kurths}}
\affiliation{Potsdam Institute for Climate Impact Research (PIK), Member of the Leibniz Association, P.O. Box 60 12 03, D-14412 Potsdam, Germany}
\affiliation{Department of Physics, Humboldt University of Berlin, Newtonstr. 15, 12489 Berlin, Germany}

\date{\today}

\begin{abstract}
  Time lags occur in a vast range of real-world dynamical systems due to finite reaction times or propagation speeds. Here we derive an analytical approach to determine the asymptotic stability of synchronous states in networks of coupled inertial oscillators with constant delay. Building on the master stability formalism, our technique provides necessary and sufficient delay master stability conditions. 
  We apply it to two classes of potential future power grids, where processing delays in control dynamics will likely pose a challenge as renewable energies proliferate. Distinguishing between phase and frequency delay, our method offers an insight into how bifurcation points depend on the network topology of these system designs.
\end{abstract}

\maketitle

\paragraph*{Introduction.} The study of nonlinear dynamics on complex networks has received full interdisciplinary attention in past years due to its potential for modeling the complexity of real-world dynamical systems. An intrinsic feature of such systems is that their time evolution generally depends on past states. Time delays, caused by finite propagation speeds or processing times, induce retarded reactions of variables to changes in the system. For example, delays occur in lasers because of the finite speed of light \cite{soriano-laser, ruschel-laser}; population dynamics depend on maturation and gestation times \cite{kuang-dde}, and the exchange of information between neurons requires time for both signal transmission as well as processing \cite{erneux-intro}.

Mathematically, continuous delay problems are described by delay differential equations (DDEs) \cite{bellman-dde, diekmann-dde}. From their analysis it is known that delays can substantially alter a system's asymptotic behavior \cite{olgac2002}. However, asymptotic stability analysis of DDEs is challenging because the corresponding spectrum contains an infinite number of complex roots. In fact, exact conditions for stability pose an open problem in research, especially regarding networks. Most previous studies have been limited to numerical investigations of characteristic equations or restricted to simple network topologies, often yielding only sufficient stability criteria.

Recent work has led to a thorough analytical understanding of the spectrum in the limit of large delay, with applications in e.g. optoelectronics \cite{maia-sync, sieber-large, lichtner-large, ruschel-large}.
In many cases, however, time lags may match the system's dynamical timescales and may play a critical role for stability. Particularly in systems of coupled oscillators like the paradigmatic Kuramoto model \cite{kuramoto-self}, delays often become comparable to the oscillation period. There,
the asymptotic stability of a synchronous regime is a central property with crucial implications for applications.

Pecora and Carroll have developed a powerful method known as the \textit{master stability formalism} to determine the stability of synchronization for identical oscillators without delay \cite{pecora-msf}. The main idea is to project the state vector into the eigenspace of the coupling matrix, yielding a block diagonal form that defines the associated master stability function (MSF). This way, dynamical parameters of the system are separated from topological information about the network.

Several studies have calculated MSFs for specific models with time-delayed couplings \cite{dahms-dmsf, dhamala-dmsf, kinzel-dmsf, lehnert-dmsf}. Here, for the first time, we generalize the formalism to DDE inertial oscillator models containing an arbitrary constant discrete delay $\tau > 0$ that may appear in the local dynamics as well as in a diffusive coupling term. While the master stability formalism requires complete synchronization of oscillators, we merely assume phase synchronization where oscillators may have constant phase differences \cite{pikovski-sync}. Our analytic approach leads to a decomposition into second-order DDEs in terms of the eigenvalues of the graph Laplacian matrix. Based on results from Bhatt and
Hsu \cite{bh}, we derive necessary and sufficient conditions for the asymptotic stability of synchronized inertial oscillator networks with delay. The corresponding \textit{delay master stability function} (dMSF) is given in terms of the graph Laplacian spectrum, the delay $\tau$, as well as dynamical parameters of the model.

For delays caused by processing times, our results offer a complete analytic solution to the question of asymptotic stability. The dMSF comprises a finite number of easily evaluated critical conditions that hold for any network topology. Particularly, in an important case which covers our central application of renewable inverter-based power grids, the conditions further simplify to a single stability criterion involving just the maximum graph Laplacian eigenvalue.

\paragraph*{Main application.} We begin with an inverter-based power grid model to exemplify how we obtain a concise condition for a major application of oscillator networks. Due to the energy transition, power grids currently undergo substantial structural and dynamical changes, threatening stable synchronization of the AC voltage frequency \cite{kroposki-pg, rohden2012self}. Characterized by a large share of volatile distributed generation units, e.g. solar or wind power plants, future energy networks will require novel control approaches like grid-forming power inverters to maintain stability \cite{schiffer-survey, anvari-fluctuations}. As this involves measurements and processing, delays are expected to play a critical role \cite{boettcher-delay, efimov2016, nussbaumer2008}. Understanding their influence on stability is thus vital to ensure security of supply and prevent blackouts.

Specifically, we consider frequency dynamics in a droop-controlled inverter grid \cite{schiffer-cond} where the
steady-state power flow between two nodes depends on the sine of their phase difference,
\begin{align} \label{eq:inverter-model}
  \ddot{\varphi}_i &= - \tilde \alpha \dot \varphi_i + \tilde \beta \bigg( P_i^d - \sum_{j=1}^N K_{ij} \ \sin \big( \Delta \varphi_{ji}^\tau  \big) \bigg) \ .
\end{align}
Here $\varphi_i(t)$ denotes the phase angle of the $i$-th inverter (oscillator) and $\Delta \varphi_{ji}^\tau := \varphi_i(t-\tau) - \varphi_j(t-\tau)$. $\tilde \alpha > 0$ and $\tilde \beta > 0$ are the inertia-specific damping and droop constants, respectively; $P_i^d$ represents the desired active power set points. Elements of the weighted adjacency matrix $(K_{ij})$ may be interpreted as the maximally transmittable power values along transmission lines in the network \cite{machowski-pg} (details in SI \footnote{See Supplemental Material (attached below) for technical details, an elaborate derivation, and more on decisive roots. Further information is also available at \url{https://github.com/reykboerner/delay-networks}.}).

We find that a synchronous state of Eq. \eqref{eq:inverter-model} is asymptotically stable if and only if
\begin{align} \label{eq:inverter-msc-final}
  \lambda_N < \frac{1}{\tilde{\beta}} \, \sqrt{\left(\frac{y_1}{\tau} \right)^4  + \tilde{\alpha}^2 \left(\frac{y_1}{\tau} \right)^2} \ , \quad y_1 = \tilde \alpha \tau \cot y_1 \ ,
\end{align}
where $y_1 \in (0,\pi]$. This exact stability condition depends on network structure only via the largest eigenvalue $\lambda_N$ of the effective Laplacian matrix $\mathcal{L}$ (see below). Notably, it suffices to compute precisely one unique characteristic root $y_1$ of the linearized spectrum associated with Eq. \eqref{eq:inverter-model}. We discuss this result further after deriving the general approach.

\paragraph*{Derivation.}
We consider a nonlinear dynamical system of $N$ coupled oscillators on a network. All oscillators (nodes) have inertia, obeying a Newtonian law of motion. The state of the $i$-th oscillator at time $t$ is given by the phase angle $\varphi_i(t)$ and angular frequency deviation $\omega_i(t) \equiv \Dot{\varphi}_i(t)$ in a reference frame co-rotating with a coherent frequency $\Omega$. Let the time evolution of the global state be governed by a set of second-order DDEs containing a discrete, constant delay $\tau > 0$,
\begin{align} \label{eq:model}
  \ddot{\varphi}_i =& \ f_i(\varphi_i, \dot \varphi_i) + f^\tau(\varphi_i^\tau, \dot \varphi_i^\tau) \nonumber \\
  & \ + \sum_{j=1}^N A_{ij} \big[ g(\Delta \varphi_{ij}, \Delta \dot \varphi_{ij}) + g^\tau(\Delta \varphi_{ij}^\tau, \Delta \dot \varphi_{ij}^\tau) \big] \ ,
\end{align}
for $i \in \{ 1, \dots, N \}$. Here time arguments are abbreviated as $\varphi_i \equiv \varphi_i(t)$ and $\varphi_i^\tau \equiv \varphi_i(t-\tau)$; furthermore $\Delta \varphi_{ij} \equiv \varphi_j - \varphi_i$ and $\Delta \varphi_{ij}^\tau \equiv \varphi_j^\tau - \varphi_i^\tau$.  The real scalar functions $f_i$ and $f^\tau$ represent undelayed and delayed isolated dynamics, respectively. Unlike identical oscillators, $f_i$ may differ from node to node by an additional constant $c_i \in \mathbb{R}$, which accounts for heterogeneous driving forces. In the interaction term, $g$ denotes undelayed coupling dynamics, whereas $g^\tau$ describes interactions with a coupling processing delay. The strength of the coupling as well as the network topology are stored in the weighted adjacency matrix $A \in \mathbb{R}^{N \times N}$, with $A_{ij} > 0$ if nodes $i$ and $j$ are connected and 0 otherwise. Here, we consider undirected graphs without self-loops.

The delay considered in Eq. \eqref{eq:model} is a processing delay which arises, for example, in engineered systems with feedback control due to measurement and processing times. Contrarily, transmission or communication delays in diffusive coupling of the form $g^{\tau}(x_j^\tau-x_i)$ require separate treatment (see SI).

To assess asymptotic stability, we linearize our DDE model near the phase synchronization manifold $\mathcal{Z}$, defined by
$\mathcal{Z} := \{ (\varphi_i, \omega_i) \in \mathbb{R}^2 : \omega_i = \dot \omega_i = 0 \ \forall i \}$.
Physically, this means that all oscillators are entrained to a coherent frequency $\Omega$ but possibly with fixed phase differences between them. A synchronous solution $\varphi^* = (\varphi_1^*, \dots, \varphi_N^*)$ with $\omega^* = (0,\dots,0)$ lies on $\mathcal{Z}$ and corresponds to a fixed point of Eq. \eqref{eq:model}.

Traditional MSFs require \textit{complete synchronization}, i.e. all oscillators move in phase with frequency $\Omega$. Then, Jacobians evaluated on the synchronization manifold are identical for all nodes \cite{pecora-msf}. To achieve a block decomposition similar to MSF for \textit{phase synchronization}, we assume:
\textit{1)} The Jacobian matrices of all local functions $f_i$ and $f^\tau$, evaluated on $\mathcal{Z}$, are identical.
\textit{2)} The Jacobians of all coupling functions $g$ and $g^{\tau}$, evaluated on $\mathcal{Z}$, are edge-independent except for a pre-factor $\text{w}_{ij}(\Delta \varphi_{ij}^*) = \text{w}_{ji}(\Delta \varphi_{ji}^*) \in \mathbb{R}$ which may depend on the fixed point $\varphi^*$.
We note that the following procedure also holds for more general coupling functions $g(\varphi_i, \dot \varphi_i, \varphi_j, \dot \varphi_j)$ and $g^{\tau}(\varphi_i^\tau, \dot \varphi_i^\tau, \varphi_j^\tau, \dot \varphi_j^\tau)$ if their first partial derivatives are antisymmetric with respect to the exchange of $i$ and $j$.
Nonetheless, we present the widely applied diffusive form here and refer to the SI for more information.

We define the effective Laplacian matrix $\mathcal{L}$ of the linearized network model such that $\mathcal{L}_{ij} := - \text{w}_{ij} A_{ij} + \delta_{ij} \sum_j \text{w}_{ij} A_{ij}$. This matrix is symmetric, positive-semidefinite, and consequently diagonalizable. In the spirit of MSF, we now transform coordinates into the space spanned by the eigenvectors of $\mathcal{L}$, with corresponding eigenvalues $\lambda_k$, $k \in \{ 1, \dots, N \}$.
Diagonalization does not affect the Jacobians (which are node-independent by assumption after absorbing $\text{w}_{ij}$ in the adjacency matrix), such that the system of DDEs decomposes into $N$ blocks given in terms of $\lambda_k$,
\begin{align} \label{eq:2dde}
  \Ddot{\theta}_k &= - a_k(\lambda_k) \dot{\theta}_k - b_k(\lambda_k) \theta_k - a_k^\tau(\lambda_k) \dot{\theta}_k^\tau - b_k^\tau(\lambda_k) \theta_k^\tau \ .
\end{align}
Here the set of $\theta_k$ denotes (small) phase angle deviations from $\varphi^*$ expressed in the transformed coordinates. The coefficients are given by elements of the Jacobian matrices (see SI); in the following we suppress their dependence on $\lambda_k$.

The stability of a synchronous state $\varphi^*$ depends on the real parts of the roots of the characteristic equation associated with Eq. \eqref{eq:2dde},
\begin{align} \label{eq:H}
  \mathcal{H}(z) := (z^2 + a_k \tau z + b_k \tau^2) e^z + a_k^\tau \tau z + b_k^\tau \tau^2 = 0 \ .
\end{align}
The exponential polynomial $\mathcal{H}$ features an infinite number of complex roots; all must have negative real parts for asymptotic stability. We now assume that the delay $\tau$ appears either in the time argument of the phases $\varphi_i$ or of the frequencies $\omega_i$ but not in both. For these cases, Bhatt and Hsu \cite{bh} derive necessary and sufficient stability conditions for scalar second-order DDEs, determined by a finite number of \textit{decisive roots} within the infinite spectrum of $\mathcal{H}$ (see also \cite{pontrjagin-1, pontrjagin-2}). After the block decomposition outlined above, we may transfer these conditions to inertial oscillator networks to obtain delay master stability conditions.

If only the coupling is delayed ($f^\tau = 0$), the longitudinal eigenvalue $\lambda_1 = 0$ describes dynamics within the synchronization manifold and asymptotic stability is determined by the $N-1$ transversal directions \cite{pecora-msf}. Contrarily, if $f^\tau \neq 0$, all $k$ must be considered for stability analysis. We define the transversal set $\mathcal{N}$, which is $\{ 2, \dots, N \}$ for $f^\tau = 0$ and $\{ 1, \dots, N \}$ otherwise.

Substituting $z = iy$ in Eq. \eqref{eq:H}, $\mathcal{H}$ separates into a real and an imaginary part. First, we consider the case of \textit{phase delay} ($a_k^\tau = 0$). Let $a_k > 0$ and $-a_k < b_k\tau \leq 0$. Then, for each $k$, there exists one decisive root $y_{1,k} \in (0,\pi]$ which solves the imaginary part of Eq. \eqref{eq:H} \cite{bh}. A synchronous fixed point of Eq. \eqref{eq:model} with phase delay is asymptotically stable if and only if, for all $k \in \mathcal{N}$,
  \begin{align} \label{eq:msc-1}
    &- b_k  <  b_k^\tau  < \frac{R_k(y_{1,k})}{\tau^2} \ ,
  \end{align}
with $R_k(y) := \sqrt{(y^2 - b_k\tau^2)^2 + (a_k \tau y)^2}$.

In the \textit{frequency delay} case ($b_k^\tau = 0$), we assume $a_k > 0$ and $b_k > 0$. Here we examine positive solutions $y_{k}$ of the real part of Eq. \eqref{eq:H}. The first positive root $y_{k,0}$ lies in the interval $(0,\pi/2)$, and one root $y_{k,m}$ is situated in each $\pi$-interval $(m\pi-\pi/2, m\pi+\pi/2)$ for $m = 1,2,\dots$. Of these roots, the decisive roots $y_k^*$ and $y_k^{**}$ are found according to $y_k^* = \underset{m \text{ odd}}{\min} \big| y_{k,m} - \tau \sqrt{b_k} \big| \text{ and } \underset{m \text{ even}}{\min} \big| y_{k,m} - \tau \sqrt{b_k} \big|$. Then, it is necessary and sufficient for asymptotic stability of a synchronous state that, for all $k \in \mathcal{N}$,
\begin{align} \label{eq:msc-2}
  &- \frac{R_k(y_k^{**})}{y_k^{**}} < a_{k}^\tau \tau < \frac{R_k(y_k^{*})}{y_k^{*}} \ ,
\end{align}
where $R_k(y)$ is defined beneath Eq. \eqref{eq:msc-1}.
Though finding decisive roots appears more complicated than in the previous case, it turns out that we must find at maximum $2N$ roots in total, within known intervals. Details are provided in the SI.

\paragraph*{Phase delay.} Our motivating example introduced in Eq. \eqref{eq:inverter-model} illustrates the case of phase delay. Here the root $y_{1,k}$ is $k$-independent and we must only consider the largest eigenvalue $\lambda_N$ of $\mathcal{L}$. We may write the resulting stability condition (Eq. \eqref{eq:inverter-msc-final}) as a dMSF $\sigma$,
\begin{align}
    \sigma(\lambda_N, \tau) = \lambda_N - \frac{1}{\tilde{\beta}} \, \sqrt{\left(\frac{y_1}{\tau} \right)^4  + \tilde{\alpha}^2 \left(\frac{y_1}{\tau} \right)^2} \ .
\end{align}
A combination of $\lambda_N$ and $\tau$ is stable if and only if $\sigma < 0$. Monotonicity arguments for $y_1(\tau)$ prove the existence of precisely one critical delay $\tau_c = \tau_c(\tilde \alpha, \tilde \beta, \lambda_N)$. A synchronous state that is stable without delay will remain asymptotically stable for all $\tau < \tau_c$ and is unstable for all $\tau \geq \tau_c$.
\begin{figure}[b]
  \centering
  \includegraphics[width=\columnwidth]{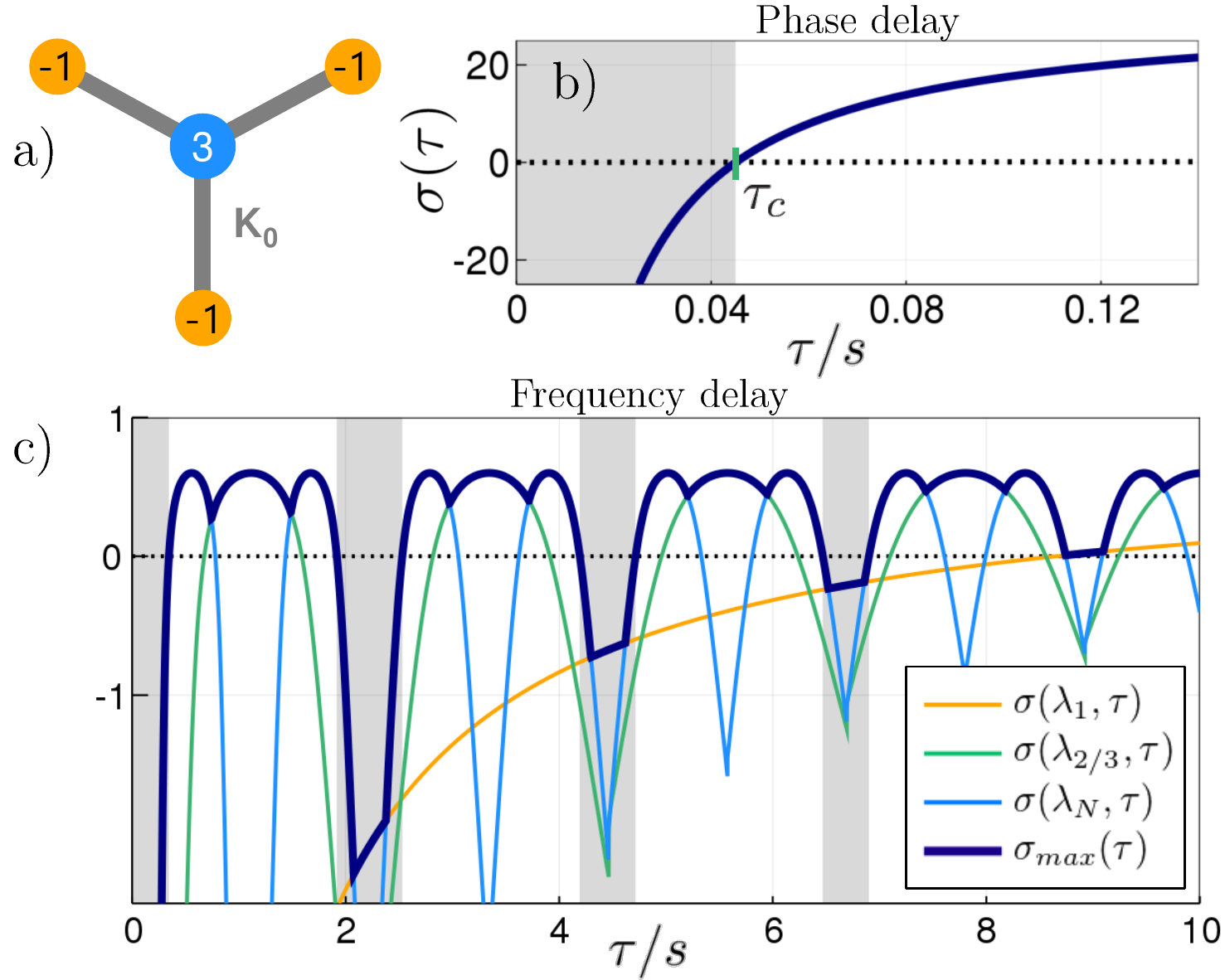}
  \caption{(Color online). Delay master stability functions for a four-node star topology (a) with three consumers connected to one producer. (b) dMSF as a function of delay $\tau$ for the inverter model with phase delay. The system is asymptotically stable for all $\tau<\tau_c \approx 45$ ms (grey area). (c) dMSF for the DSGC model with frequency delay. Each Laplacian eigenvalue $\lambda_k$ contributes a curve $\sigma(\lambda_k,\tau)$; the system is stable in regimes where $\sigma_{max}<0$ (grey areas). Parameter values are $P_0 = 1 \hzz$, $K_0=8 \hzz$, $\tilde{\alpha}=\alpha =0.1 \hz$, $\tilde{\beta}=0.07$, and $\gamma = 0.25 \hz$.}
  \label{fig:msf}
\end{figure}

\begin{figure}
  \centering
  \includegraphics[width=\columnwidth]{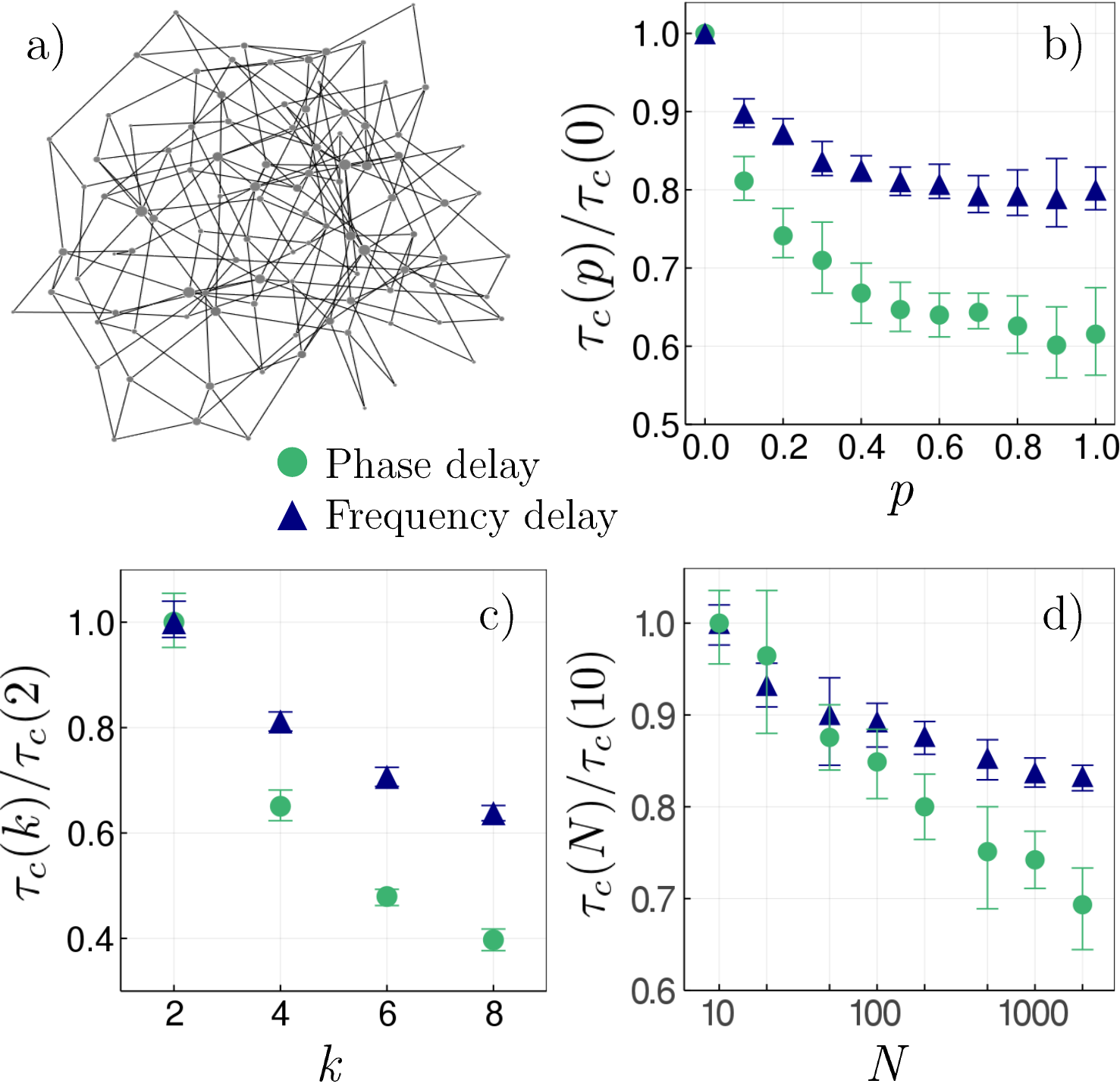}
  \caption{(Color online). Critical delay $\tau_c$ on Watts-Strogatz networks. (a) shows a network example for $N=100$ nodes, degree $k=4$, and rewiring proabability $p=0.5$. Consumers/producers with equal power input/output, connected via lines of equal admittance, are placed alternately on the original ring. Plots (b-d) illustrate how $\tau_c$ depends on $p,k,N$ for the inverter model (phase delay, green circles) and for the DSGC model (blue triangles). Parameter values not explicitly given are $N=100$, $k=4$, $p=0.5$; model parameters are as in Fig. 1. Vertical axes are normalized for each series; data points are averaged over 10 realizations.}
  \label{fig:ws}
\end{figure}

To visualize $\sigma$, we first choose a star topology as in Ref. \cite{schaefer-taming}. A producer in the center with steady-state power production $P_+ = 3P_0$ is connected to three consumers with $P_- = -P_0$ via transmission lines of equal capacity $K_0$ (Fig. \ref{fig:msf}a). Due to the symmetry of the configuration, we obtain three distinct Laplacian eigenvalues. The phase delay case depends on $\lambda_N$ only; thus we get a single curve with $\sigma=0$ at $\tau_c$ (Fig. \ref{fig:msf}b). For typical parameter values, $\tau_c \approx 40$ ms is about twice the 50 Hz oscillation period.

\paragraph*{Frequency delay.} Stability in the presence of a frequency delay is qualitatively different. We show this by applying our approach to the decentral smart grid control (DSGC) scheme \cite{schaefer-decentral, schaefer-taming}. The model incorporates electricity price dynamics by locally relating the price to the current grid frequency, motivating producers/consumers to adapt their feed-in/consumption to the currently available power supply. The continuous measurements required for this smart grid regulation induce a local processing delay at each node. The model reads \cite{schaefer-taming}
\begin{equation} \label{schaefer-model}
      \Ddot{\varphi}_i = P_i  - \alpha \Dot{\varphi}_i - \gamma \Dot{\varphi}_i^\tau + \sum_j K_{ij} \sin(\varphi_j - \varphi_i)  \ ,
\end{equation}
where $P_i$ is the produced/consumed power at node $i$ and $K = (K_{ij})$ denotes the weighted adjacency matrix as before. Next to the damping $\alpha > 0$, the price elasticity $\gamma > 0$ acts as a second, delayed damping term.

For the DSGC model we obtain $N$ delay master stability conditions
\begin{align} \label{eq:dmsf-2}
  0 < \gamma \tau < \frac{1}{y_k^*} \sqrt{({y_k^*}^2 - \lambda_k \tau^2)^2 + (\alpha \tau y_k^*)^2} \ .
\end{align}
Here it is not a priori identifiable which root $y_k^* = y_k^*(\lambda_k, \tau)$ determines stability for a given $\tau$; we must regard all Laplacian eigenvalues $\lambda_k$. The root $y_k^{**}$ is not relevant because $\gamma \tau > 0$. Analogous to the previous case, we may formulate Eq. \eqref{eq:dmsf-2} as a dMSF $\sigma(\lambda_k,\tau)$. Then, the system is stable in all regions where $\sigma_{max}(\tau) := \underset{\lambda_k}{\max} \  \sigma(\lambda_k, \tau) < 0$.

Calculating $\sigma_{max}$ for the star topology (Fig. \ref{fig:msf}a), we now have contributions from all three distinct eigenvalues $\lambda_k$ as depicted in Fig. \ref{fig:msf}c. In addition to a stable regime beginning at $\tau=0$, there exist further windows of stability for larger delays, corroborating prior results based on numerical analysis \cite{schaefer-taming}.

\paragraph*{Networks.}  For frequency delays, we conjecture that the length of the first stability window, extending from $\tau=0$ to a critical delay $\tau_c$, is determined by $\lambda_N$. In the limit of an infinite, heterogenous graph, its Laplacian spectrum may be expected to become quasi-continuous, such that further stability windows vanish. Thus, in both frequency and phase delay, the maximum eigenvalue of the effective graph Laplacian plays a crucial role. Several bounds and estimates in terms of network characteristics have been published for the largest Laplacian eigenvalue of a weighted graph, e.g. \cite{sorgun-bounds}. Therefore, our method may provide insight even without explicitly diagonalizing $\mathcal{L}$, which could be practical particularly for large systems.

Finally, we explore how the critical delay $\tau_c$ (end of first stability window) depends on the network structure in both models. As a versatile example, we generate Watts-Strogatz networks \cite{Watts1998} with different rewiring probabilities $p$, mean degrees $k$, and number of nodes $N$ (e.g. Fig. \ref{fig:ws}a). Consumers and producers are placed alternately on the original ring graph. The results for varying $p,k,$ and $N$ (Figs. \ref{fig:ws}b-d) show that the critical delay decreases with an increasing number of nodes and edges as well as randomness of the system, which emphasizes the significance of delays for the design and control of real-world dynamical systems. However, the decline is smaller for the DSGC model compared to the droop-controlled inverter model.

\paragraph*{Conclusion.}
In this Letter, we present an analytical approach to assess the asymptotic stability of synchronized inertial oscillator networks with delayed dynamics. Specifically, we consider processing delays either in the phase or in the frequency. We show how to extend the master stability formalism to an arbitrary lag time, obtaining dMSFs in terms of the delay $\tau$ and eigenvalues of the effective graph Laplacian $\mathcal{L}$. Unlike MSF, we more generally consider phase synchronization and allow for constant inhomogeneities in the local dynamics of otherwise identical oscillators. A block decomposition of the linearized model yields necessary and sufficient stability conditions which deliver an analytic expression for the dependence on network structure. These criteria involve maximally $2N$ decisive roots. In contrast, previous numerical asymptotic stability analyses rely on randomly computing a significantly larger number of characteristic roots without being certain that all decisive roots have been found. Illustrating our approach, we consider two concrete models for renewable power grids as our main applications. Notably, we are able to boil down the problem of stability to a single condition in the case of the droop-controlled inverter model. We therefore believe that our method could contribute to the development of design criteria for future energy systems. Generally, our results advance the stability analysis of dynamical systems on complex oscillator networks complying with Eq. \eqref{eq:model}.

Due to its increasing importance in real-world applications, delay stability in complex systems remains an important topic with many open challenges. Our work opens a new analytic approach on this subject. From the perspective of power grids, it is crucial to also tackle non-identical oscillators, to consider multiple delays, and combined phase and frequency delays. In the wider context of physical systems, it is also highly interesting to study more general stability criteria for diffusive coupling with transmission delays.

The authors acknowledge the support of BMBF, Condynet2 FK. 03EK3055A.
This work was funded by the Deutsche Forschungsgemeinschaft (DFG, German Research Foundation)
– KU 837/39-1 / RA 516/13-1. All authors gratefully acknowledge
the European Regional Development Fund (ERDF), the German Federal Ministry of
Education and Research and the Land Brandenburg for supporting this project by
providing resources on the high performance computer system at the Potsdam
Institute for Climate Impact Research.

\bibliography{boerner_dmsf}

\onecolumngrid
\newpage

\renewcommand\theequation{S.\arabic{equation}}
\renewcommand\thefigure{S.\arabic{figure}}
\setcounter{equation}{0}
\setcounter{figure}{0}

\begin{center}
    \Large SUPPLEMENTAL INFORMATION (SI)
\end{center}

\section{Coefficients of second-order DDE blocks}
The coefficients in Eq. (4) of the main text are given by
\begin{align} \label{eq:second-order-block}
  a_k &= F_\omega - \lambda_k G_\omega \  &\quad a_k^\tau &= F_\omega^\tau - \lambda_k G_\omega^\tau \nonumber \\
  b_k &= F_\varphi - \lambda_k G_\varphi \  &\quad b_k^\tau &= F_\varphi^\tau - \lambda_k G^{\tau}_\varphi \ ,
\end{align}
where
\begin{align}
    F_{\varphi} &:= \eval{\pdv{f_i(\varphi_i, \dot \varphi_i)}{\varphi_i}}_{\varphi^*} \ &\quad  G_{\varphi} &:= \frac{1}{\text{w}_{ij}(\varphi^*)} \eval{\pdv{g(\Delta \varphi_{ij}, \Delta \dot \varphi_{ij})}{(\Delta \varphi_{ij})}}_{\varphi^*} \nonumber \\
    F_{\omega} &:= \eval{\pdv{f_i(\varphi_i, \dot \varphi_i)}{\dot \varphi_i}}_{\varphi^*} \ &\quad G_{\omega} &:= \frac{1}{\text{w}_{ij}(\varphi^*)} \eval{\pdv{g(\Delta \varphi_{ij}, \Delta \dot \varphi_{ij})}{(\Delta \dot \varphi_{ij})}}_{\varphi^*} \ .
\end{align}
Elements of the delayed Jacobians are written analogously as $F_\varphi^\tau$ etc.

\section{General Derivation}
In the Letter, we outline the derivation of our approach based on the inertial oscillator model described by Eq. (3) of the main text. Here we present a more elaborate variant (discussed in detail in Ref. \cite{boerner}). This allows us to illuminate the underlying assumptions of our method and discuss why including communication delays involves a strong restriction.

Recall that we formulate the oscillators' dynamics in a reference frame co-rotating with the frequency $\Omega$. A synchronous state where all oscillators are entrained to frequency $\Omega$ corresponds to a fixed point characterized by $\omega_i = \dot \omega_i = 0 \ \forall i$.

Instead of the second-order form stated in Eq. (3) of the main text, we may equivalently express our inertial oscillator network model as a set of first-order DDEs by treating the phase angles $\varphi_i(t)$ and angular frequency deviations $\omega_i(t) = \dot \varphi_i(t)$ as two independent variables for each node. We thus define the vector $x_i \equiv (\varphi_i, \omega_i)^\top$ and write
\begin{align}\label{model-delay}
    \Dot{x}_i = \ & \boldsymbol f_i(x_i) + \boldsymbol f^\tau(x_i^\tau) + \sum_{j=1}^N A_{ij}  \left( \boldsymbol g^{00}(x_i, x_j) + \boldsymbol g^{\tau\tau}(x_i^\tau, x_j^\tau) + \boldsymbol g^{0\tau}(x_i, x_j^\tau) \right) \ ,
\end{align}
where we have abbreviated $x_i \equiv x_i(t)$ and $x_i^\tau \equiv x_i(t-\tau)$. Here we denote vector-valued functions ($\mathbb{R}^2 \to \mathbb{R}^2$) by bold letters, while $f_i,f^\tau,g$, and $g^\tau$ will remain the scalar functions introduced in the main text; e.g. $\boldsymbol f_i(x_i) = (\omega_i,f_i(\varphi_i,\omega_i))^\top$. In contrast to Eq. (3), this form now includes a \textit{communication delay} via the function $\boldsymbol g^{0\tau}(x_i(t),x_j(t-\tau))$. Furthermore, note that the coupling functions $\boldsymbol g^{00}, \boldsymbol g^{\tau \tau}$, and $\boldsymbol g^{0\tau}$ may depend on $x_i$ and $x_j$ in an arbitrary fashion.

A communication delay of the type above may describe transmission or propagation lags between nodes. The intuition is that a change of node $i$ at time $t$ depends on the history of connected nodes because it takes the time $\tau$ until a signal from a node $j$ reaches node $i$.

We now linearize Eq. \eqref{model-delay} around a fixed point $x^* = (x_1^*, \dots, x_N^*)^\top$ defined by the conditions $\dot{x}_i^* = (0,0)$ for all $i$. The set of fixed points constitutes the phase synchronization manifold $\mathcal{Z}$,
$$ \mathcal{Z} := \{ x_i \in Z \subset \mathbb{R}^2 : \quad i = 1, \dots, N \ \text{and} \ \dot x_i = (0,0) \ \forall i \} \ .$$
Substituting $x_i = x_i^\ast + \eta_i$, with $|| \eta_i ||$ small for all $i$, this yields
\begin{align} \label{model-linear}
    \Dot{\eta}_i \approx & \ \text{D}_i \boldsymbol  f_i \eta_i + \text{D}_i \boldsymbol f^\tau \eta_i^\tau \\
    &+ \sum_{j=1}^N A_{ij} \left[
    \left( \text{D}_{ij}^1 \boldsymbol g^{00} \eta_i + \text{D}_{ij}^2 \boldsymbol g^{00} \eta_j \right) +
    \left( \text{D}_{ij}^1 \boldsymbol g^{\tau \tau} \eta_i^\tau + \text{D}_{ij}^2 \boldsymbol g^{\tau \tau } \eta_j^\tau \right) +
    \left( \text{D}_{ij}^1 \boldsymbol g^{0\tau} \eta_i + \text{D}_{ij}^2 \boldsymbol g^{0\tau} \eta_j^\tau \right)
    \right] \ , \nonumber
\end{align}
where Jacobian matrices, all evaluated at the fixed point, are written in short notation,
$$ \text{D}_i \boldsymbol f := \at{\frac{\partial \boldsymbol f(x_i)}{\partial x_i}}{x_i = x_i^*} \ , \qquad  \text{D}_{ij}^{1/2} \boldsymbol g := \at{\frac{\partial \boldsymbol g(x_i, x_j)}{\partial x_{i/j}}}{x_i = x_i^*, \ x_j = x_j^*} $$
and $\eta_i^\tau \equiv \eta_i(t-\tau)$.
Due to the relation $\omega_i = \dot{\varphi}_i$ between coordinates of the vector $x_i$, some elements of the Jacobians are immediately zero or one. Particularly,
\begin{align} \label{eq:jacobians-isolated}
  \text{D}_i \boldsymbol f_i =
  \begin{bmatrix}
    0 & 1 \\
    \partial_{\varphi_i}  f_i & \partial_{\omega_i}  f_i
  \end{bmatrix}
   \ , \qquad \text{D}_i \boldsymbol f^\tau =
   \begin{bmatrix}
     0 & 0 \\
     \partial_{\varphi_i} f^\tau & \partial_{\omega_i} f^\tau
   \end{bmatrix}
   \ ,
\end{align}
where $\partial_{m} f$ denotes the partial derivative of the function $f$ by the argument $m$, evaluated at the fixed point. In the same manner,
\begin{align} \label{eq:jacobians-coupling}
  \text{D}_{ij}^1 \boldsymbol g^{00} =
  \begin{bmatrix}
    0 & 0 \\
    \partial_{\varphi_i} g & \partial_{\omega_i} g
  \end{bmatrix}
   \ , \qquad \text{D}_{ij}^2 \boldsymbol g^{00} =
   \begin{bmatrix}
     0 & 0 \\
     \partial_{\varphi_j} g & \partial_{\omega_j} g
   \end{bmatrix}
   \ , \qquad \text{etc.}
\end{align}
We emphasize that these Jacobians depend on the fixed point. For the phase synchronization manifold, this implies that the Jacobians may differ for different $i$ and $j$.
\subsection{Antisymmetric coupling}
Assume now that
\begin{enumerate}[start=1]
  \item there is no communication delay (i.e. $\boldsymbol g^{0\tau} = 0$) and
  \item the linearized coupling between two nodes $i$ and $j$ is antisymmetric, that is, $\text{D}_{ij}^1 \boldsymbol g^{00} (x_i^*,x_j^*) = - \text{D}_{ij}^2 \boldsymbol g^{00} (x_i^*,x_j^*)$ and $\text{D}_{ij}^1 \boldsymbol g^{\tau \tau} (x_i^*,x_j^*) = - \text{D}_{ij}^2 \boldsymbol g^{\tau \tau} (x_i^*,x_j^*)$.
\end{enumerate}
This is fulfilled by the model discussed in the main text (Eq. (3)), where we have diffusive coupling. The antisymmetry requirement will allow us to write the problem in terms of the effective graph Laplacian matrix $\mathcal{L}$.

Our goal is to decouple local information about the dynamics of single nodes from global terms characterizing the network as a whole. Mathematically, this is achieved when local $2\times 2$ matrices and global $N \times N$ matrices factorize into a Kronecker product (symbolized by $\otimes$).

If we have complete synchronization (i.e. all nodes oscillate with identical frequency \textit{and} phase angle), the Jacobians are homogeneous for all $i,j$, resulting in immediate Kronecker factorization. This is not true for the more general case of phase synchronization. To achieve a decomposition nonetheless, we impose the following restrictions in analogy to the main text:
\begin{enumerate}[start=3]
    \item The Jacobians of the local functions $f$ and $f^\tau$, respectively, evaluated on the phase synchronization manifold, are identical for all nodes: $$ F : = \text{D}_1 \boldsymbol f = \text{D}_2 \boldsymbol f = \cdots = \text{D}_N \boldsymbol f \  $$
    $$ F^\tau : = \text{D}_1 \boldsymbol f^\tau = \text{D}_2 \boldsymbol f^\tau = \cdots = \text{D}_N \boldsymbol f^\tau \ .$$
    \item The Jacobians of the coupling functions $\boldsymbol g^{00}$ and $\boldsymbol g^{\tau\tau}$, respectively, evaluated on the phase synchronization manifold, are identical for all $i,j$ up to a prefactor $\text{w}_{ij}(x_i^*,x_j^*) \in \mathbb{R}$ which contains all dependencies on the fixed point. It has the property $\text{w}_{ij} = \text{w}_{ji} > 0$;
    $$ \text{w}_{ij}(x_i^*,x_j^*) G^{00} := \text{D}^2_{ij} \boldsymbol g^{00}(x_i^*,x_j^*) $$
    $$ \text{w}_{ij}(x_i^*,x_j^*) G^{\tau\tau} := \text{D}^2_{ij} \boldsymbol g^{\tau\tau}(x_i^*,x_j^*) \ . $$
\end{enumerate}
With these assumptions, all dependencies on the fixed point may be absorbed in the \textit{effective adjacency matrix} $\mathcal{A}$ with entries $\mathcal{A}_{ij} := \text{w}_{ij} A_{ij}$. Since we assume antisymmetric coupling (assumption 2), we may furthermore replace $\mathcal{A}$ by the the \textit{effective graph Laplacian matrix} $\mathcal{L}$ given by $\mathcal{L}_{ij} := \delta_{ij} \sum_j \mathcal{A}_{ij} - \mathcal{A}_{ij}$. Now, the set of linearized DDEs reads
\begin{align}
\Dot{\eta_i} = F \eta_i + F^\tau \eta_i^\tau - \sum_{j=1}^N \mathcal{L}_{ij} \Big( G^{00} \eta_j + G^{\tau \tau} \eta_j^\tau \Big) \ ,
\end{align}
or, in vector notation for the entire system, $\eta = (\eta_1, \dots, \eta_N)^\top$,
\begin{align}
    \Dot{\eta} = [ \mathbb{I}_N \otimes F - \mathcal{L} \otimes G^{00}] \eta + [ \mathbb{I}_N \otimes F^\tau - \mathcal{L} \otimes G^{\tau\tau}] \eta^\tau \ .
\end{align}
Here $\mathbb{I}_N$ is the $N$-dimensional unit matrix. According to assumption 4, $\mathcal{L}$ is symmetric and therefore diagonalizable. Switching to a basis $\mathcal{B}$ of eigenvectors via the coordinate transform $\xi = [T_\mathcal{B} \otimes \mathbb{I}_2] \eta$, we diagonalize $\mathcal{L} = T_\mathcal{B}^{-1}\Lambda T_\mathcal{B}$ to obtain the diagonal matrix of its eigenvalues, $\Lambda = diag(\lambda_1, \dots, \lambda_N)$. This leads to a block-diagonal form; each two-dimensional block is given by the equation
\begin{align} \label{eq:blocks}
    \Dot{\xi}_k = \big( F - \lambda_k G^{00} \big) \xi_k + \big( F^\tau - \lambda_k G^{\tau \tau} \big) \xi_k^\tau \ ,
\end{align}
where $k \in 1,\dots , N$. (Note that we have switched index from $i$ to $k$ to emphasize that $\xi$ represents the state vector in the eigenbasis $\mathcal{B}$.)

Similar to the master stability formalism, we have thus decomposed the problem into blocks which vary only in the effective Laplacian eigenvalue $\lambda_k$. All eigenvalues are non-negative and $\lambda_1 = 0$ because, according to the properties of an undirected graph's Laplacian matrix, $\mathcal{L}$ is positive-semidefinite.

Recalling that the two components of the vector $x_i$ are related via $\omega_i = \dot \varphi_i$, we may write Eq. \eqref{eq:blocks} in second-order form. The vector $\phi$ of all phase angles is given by the linear combination $\phi = \sum_k v_k \theta_k$, where $v_k$ is an eigenvector of $\mathcal{L}$. In terms of the phase angles in the transformed coordinates, $\theta_k$, we obtain a second-order DDE,
\begin{align} \label{eq:2order}
  \Ddot{\theta}_k &= - a_k(\lambda_k) \dot{\theta}_k - b_k(\lambda_k) \theta_k - a_k^\tau(\lambda_k) \dot{\theta}_k^\tau - b_k^\tau(\lambda_k) \theta_k^\tau \ .
\end{align}
The coefficients are given in Eq. \eqref{eq:second-order-block} as functions of $\lambda_k$ and the Jacobians $F$, $F^\tau$, $G$, and $G^{\tau\tau}$. Following Bhatt and Hsu \cite{bh}, we state stability criteria in terms of these coefficients, as discussed in the main text.

\subsection{Communication delays}
The derivation above relies on the assumption of antisymmetric coupling between connected nodes $i$ and $j$ (assumptions 1 and 2). This permits rewriting the problem in terms of the effective Laplacian matrix $\mathcal{L}$, which leads to a decomposition in its eigenvalues.

Let us briefly consider symmetric or asymmetric coupling which does not satisfy assumptions 1 and 2. Notably, a communication delay immediately destroys the antisymmetry because the coupling function $g^{0\tau}(x_i, x_j^\tau)$ evaluates its arguments at different times. In this case, the problem cannot be expressed in graph Laplacian form. Instead, we consider an effective adjacency matrix $\mathcal{A}$ and an effective degree matrix $\mathcal{D}$. These matrices must then commute to allow a block decomposition.

In the following, let $\boldsymbol{g}$ represent any of the functions $\boldsymbol{g}^{00},\boldsymbol{g}^{\tau\tau}$, or $\boldsymbol{g}^{0\tau}$. When linearizing the inertial oscillator model, we obtain one Jacobian $\text{D}^2_{ij} \boldsymbol{g}(x_i^*,x_j^*)$ containing derivatives with respect to the $j$-th node, and a second Jacobian $\text{D}^1_{ij} \boldsymbol{g}(x_i^*,x_j^*)$ with derivatives by $x_i$. Replacing assumption 4, we assume:
\begin{enumerate}[start=5]
    \item The adjacency matrix $A$ and the Jacobian $\text{D}^2_{ij} g(x_i^*,x_j^*)$, evaluated on the phase synchronization manifold $\mathcal{Z}$, factorize into the direct product $\mathcal{A} \otimes G^{(2)}$ of an effective $N \times N$ adjacency matrix $\mathcal{A}$ and a universal $2 \times 2$ Jacobian $G^{(2)}$, such that the local matrix $G^{(2)}$ is the same for all nodes and only $\mathcal{A}$ depends on the indices $i,j$ of the network. \\
    Similarly, $A$ and $\text{D}^1_{ij} g(x_i^*,x_j^*)$ factorize into the direct product $\mathcal{A}^\prime \otimes G^{(1)}$.\\
    We require that $\mathcal{A}$ and $\mathcal{A}^\prime$ are symmetric matrices.
\end{enumerate}
Since $\mathcal{A}$ and $\mathcal{A}'$ are symmetric by assumption, they are diagonalizable. We define the effective degree matrix $\mathcal{D}$,
$$\mathcal{D}_{ij} :=
\begin{cases}
\tilde{d}_i := \sum_l \mathcal{A}^\prime_{il} & \text{if} \ \ j = i \\
0 & \text{otherwise}
\end{cases}
 \qquad i,j,l \in 1, \dots, N \ .
$$
Now we assume:
\begin{enumerate}[start=6]
    \item The matrices $\mathcal{A}$ and $\mathcal{D}$ commute, i.e.
$$ [ \mathcal{A}, \mathcal{D} ] = 0 \ . $$
\end{enumerate}
This implies that they are simultaneously diagonalizable. In that case,
the DDE (Eq. \eqref{model-linear}) decomposes into blocks in terms of the eigenvalues of $\mathcal{A}$ and $\mathcal{D}$. Then, the coefficients in Eq. \eqref{eq:2order} are functions of these eigenvalues instead of the eigenvalues $\lambda_k$ of the graph Laplacian matrix $\mathcal{L}$, and our method can be applied in analogy to the antisymmetric case.

Regular graphs with homogeneous weights present a special case where $ [ \mathcal{A}, \mathcal{D} ] = 0 $. For more complex network topologies, however, assumption 6 is generally not satisfied.

\section{Decisive roots}
In the case of phase delay, there is in principle one uniquely defined decisive root for each $k$. Thus, the stability analysis is generally based on calculating maximally $N$ decisive roots for a system of size $N$. In contrast, the frequency delay case generally requires identification of two decisive roots $y_k^*$, $y_k^{**}$ for each $k$. The value of these roots depends on the delay and on the coefficient $b_k$.

Assume the value $\sqrt{b_k \tau^2}$ is located in the interval $(j\pi - \pi/2, j\pi + \pi/2)$, where $j$ is a positive integer. If $j$ is an odd number, the closest decisive root $y_k^*$ is in the same interval. We then find the root $y_k^{**}$ either in the $\pi$-interval to the right ($j+1$) or to the left ($j-1$). To find all decisive roots for a given $k$, we must therefore calculate three roots and, among them, compare the two possible candidates for $y_k^{**}$.

\begin{figure}
  \centering
  \includegraphics[width=0.7\textwidth]{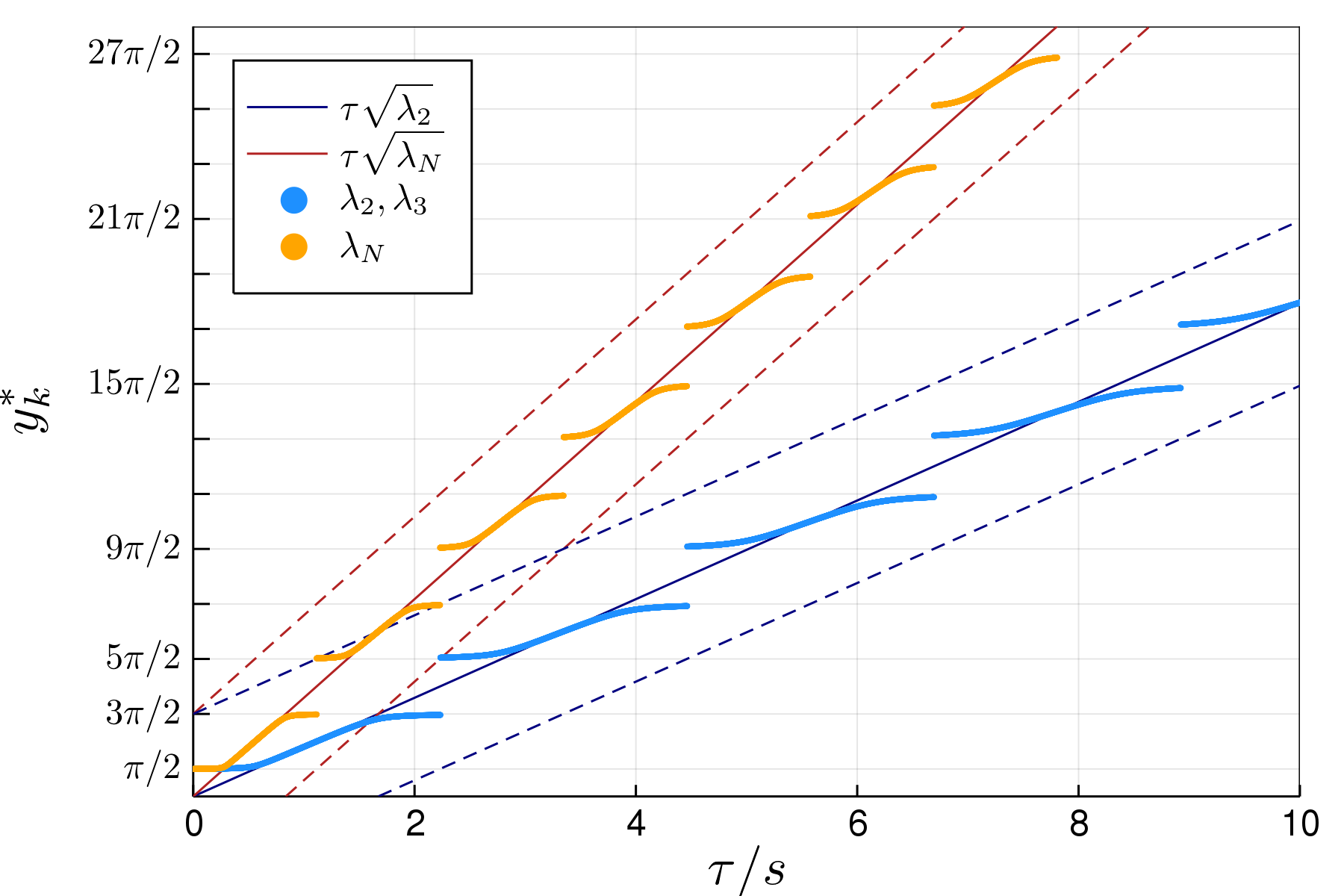}
  \caption{Scatter plot of decisive roots $y_2^*$ (blue dots) and $y_N^*$ (orange dots) as a function of delay $\tau$ for the DSGC model (frequency delay) on the four-node star network. The roots only lie in odd intervals $(j\pi-\pi/2, j\pi+\pi/2)$, $j$ odd. The solid lines show the value $\sqrt{\lambda_k \tau^2}$ for $k=2$ (blue) and $k=N$ (red). The dotted lines have a vertical distance of $3\pi/2$ from the corresponding solid line, thus bordering the interval where the root, for given $\tau$, could be found. Parameter values for the calculation are $P_0 = 1\hzz, K_0 = 8 \hzz, \alpha = 0.1 \hz$, and $\gamma = 0.25 \hz$.}
  \label{fig:roots}
\end{figure}

If $j$ is an even number, then $y_k^{**}$ is located in the same interval and $y_k^{*}$ is one of the two roots found in each of the adjacent intervals. In any case, all decisive roots must lie within a distance of $3\pi /2$ from the value $\sqrt{b_k \tau^2}$ (see figure \ref{fig:roots}).

In conclusion, for frequency delays, we must calculate at maximum $3N$ characteristic roots to find in total $2N$ decisive roots that the stability conditions are based on. In the DSGC model presented in the main text, the root $y_k^{**}$ turns out to be irrelevant; thus it suffices to calculate at maximum $2N$ characteristic roots.

\section{Droop-controlled inverter model}
In this section, we provide further detail on the renewable inverter-based power grid model with processing delay, considered in the main text as a central application of the phase delay case. For additional information on the theoretical study of power systems we refer to \cite{machowski-pg, filatrella2008analysis, dorfler-sync, hellmann-lossy}.

To model a renewable power system with phase delay, we describe the dynamics of grid-forming inverters (represented by nodes of the network) using the swing equation \cite{machowski-pg,filatrella2008analysis},
\begin{align}
    m_i \ddot \varphi_i + \alpha_i \dot \varphi_i = P_i^d - P_i^{el} \ .
\end{align}
Here $\varphi(t)$ denotes the phase angle, $m_i$ the inertia, $a_i$ the damping constant, and $P_i^d$ the desired power set point of the $i$-th node. The set point is positive for production and negative for consumption. Furthermore, we write the electrical power at node $i$ as
\begin{align} \label{eq:flow}
    P_i^{el} = \sum_j U_0^2 | B_{ij} | \sin(\varphi_i-\varphi_j) \ ,
\end{align}
where $U_0$ denotes the AC voltage amplitude which is assumed constant throughout the system, i.e. $U_0 = U_i \ \forall i$ \cite{hellmann-lossy}. $B_{ij}$ represents the susceptance of the transmission line between nodes $i$ and $j$ (we may choose its value to be zero if $i$ and $j$ are not directly connected). Eq. \eqref{eq:flow} expresses a common choice in the literature to model steady-state power flow \cite{machowski-pg}. It follows when neglecting losses (purely inductive power lines) and assuming that all phase differences $|\varphi_j-\varphi_i|<\pi/2$. 

According to Eq. \eqref{eq:flow}, the power flow along a transmission line between two nodes depends on the phase angle difference between them. Let us suppose that the state of this power line (an edge in the network) enters the frequency control, yet with a processing delay. One way to model this would be delayed coupling, where the phase difference is evaluated at time $(t-\tau)$.

Furthermore, we introduce the weighted adjacency matrix $K_{ij} := U_0^2|B_{ij}|$ and define $\tilde \alpha = \alpha_i / m_i$, $\tilde \beta = 1/m_i \ \forall i$, assuming homogeneous inertia-specific damping and droop constants. This leads to 
\begin{align}
    \ddot \varphi_i = \tilde \alpha \dot \varphi_i + \tilde \beta \bigg( P_i^d - \sum_{j=1}^N K_{ij} \ \sin \big( \Delta \varphi_{ji}^\tau  \big) \bigg) \ ,
\end{align}
which corresponds to the droop-controlled inverter model presented in Eq. (1) of the main text.

\end{document}